\documentclass{amsart}
\usepackage{amssymb,latexsym,enumerate,url}

\newtheorem{Thm}[equation]{Theorem}
\newtheorem{Prop}[equation]{Proposition}
\newtheorem{Lem}[equation]{Lemma}
\newtheorem{Cor}[equation]{Corollary}

\theoremstyle{remark}

\newtheorem*{Rem*}{Remark}
\theoremstyle{definition}

\newtheorem{Ex}[equation]{Example}
\newtheorem*{Not*}{Notation}
\newtheorem{Def}[equation]{Definition}
\numberwithin{equation}{section}

\DeclareMathOperator{\Stab}{Stab}
\DeclareMathOperator{\SL}{SL}
\DeclareMathOperator{\PSL}{PSL}
\begin{document}

\date{%
Sun May  3 23:11:01 EDT 2009}

\title{Discretely ordered groups}

\author[Peter Linnell]{Peter A. Linnell}
\address{Department of Mathematics \\
Virginia Tech \\
Blacksburg \\
VA 24061-0123 \\
USA}
\email{linnell@math.vt.edu}
\urladdr{\url{http://www.math.vt.edu/people/plinnell/}}

\author[Akbar Rhemtulla]{Akbar H. Rhemtulla}
\address{Math Dept \\
University of Alberta \\
Edmonton \\
AL Canada T6G 2G1}
\email{akbar@math.ualberta.ca}
\urladdr{\url{http://www.math.ualberta.ca/Rhemtulla_A.html}}

\author[Dale Rolfsen]{Dale P. O. Rolfsen}\
\address{Math Dept \\
University of British Columbia \\
Vancouver \\
BC Canada V6T 1Z2}
\email{rolfsen@math.ubc.ca}
\urladdr{\url{http://www.math.ubc.ca/~rolfsen/}}

\begin{abstract}
We consider group orders and right-orders which are discrete, meaning
there is a least element which is greater than the identity.  We note
that non-abelian free groups cannot be given discrete
orders, although they do have right-orders which are discrete.
More generally, we give necessary and
sufficient conditions that a given orderable group can be endowed with
a discrete order.  In particular, every orderable group $G$ embeds in
a discretely orderable group.  We also consider conditions on
right-orderable groups to be discretely right-orderable.  Finally, we
discuss a number of illustrative examples involving discrete
orderability, including the Artin braid groups and Bergman's
non-locally-indicable right orderable groups.
\end{abstract}



\maketitle

\section{Introduction}

Let $G$ be a group and suppose $<$ is a strict
total order relation on the set of its elements.  Then $(G,<)$ is a
right-ordered group if $f < g$ implies that
$fh < gh$ for all $f,g,h\in G$.
If in addition $f < g$ implies that $hf < hg$, then we say
$(G,<)$ is an ordered group.  If such an order exists for a given
group $G$, we say that $G$ is right-orderable or orderable,
respectively.  We call the order $<$ \emph{discrete} if there is an
element $a \in G$ such that $1 < a$, where $1$ denotes the identity
element of $G$, and there is no element of $G$ strictly between these.

For a right-ordered group, the positive cone $P := \{g\in G \mid 1 <
g\}$ satisfies
\begin{enumerate}
\item
$P$ is closed under multiplication and

\item
for every $g \in G$, exactly one of $g = 1$, $g \in P$ or $g^{-1}
\in P$ holds.
\end{enumerate}

Conversely, if a group $G$ has a subset $P$ with properties (1) and
(2), it is routine to verify that the order defined by $g < h$ if and
only if $hg^{-1} \in P$ makes $(G, <)$ a right-ordered group.
Similarly, a group $G$ is orderable if and only if it admits a subset
$P$ satisfying (1), (2) and

\begin{enumerate}
\item[(3)]
$gPg^{-1} = P$ for all $g \in G$.
\end{enumerate}

A subset $X$ of a (right-) ordered group $G$ is
\emph{convex} if $x < y < z$ and
$x, z \in X$ imply $y \in X$.  We recall that the set of all convex
subgroups of a (right-) ordered group is linearly ordered by
inclusion.  A convex \emph{jump} $C \rightarrowtail D$
is a pair of distinct convex
subgroups such that $C \subset D$ and there are no convex subgroups
strictly between them.  In particular, the convex jump determined by
a nonidentity element $g$ of $G$ is defined by $C =$ the union of all
convex subgroups not containing $g$ and $D =$ the intersection of all
convex subgroups which do contain $g$.
If the group is orderable, then for any convex
jump $C \rightarrowtail D$, $C$ is normal in $D$ and the quotient
$D/C$ embeds in $\mathbb{R}$, the additive reals,
by an order-preserving isomorphism.

\begin{Lem}\label{Ldense}
If $<$ is a discrete right-order on $G$ and
$a$ is the least positive element under
$<$, then $\langle a \rangle$ is convex.  Moreover, for any $g \in G$,
we have $a^{-1}g < g < ag$ and there is no
element strictly between these elements of $G$.
If the right-order $<$ is not discrete,
then it is dense in the sense that for any $f, g \in
G$ with $f < g,$ there exists $h \in G$ with $f < h < g$.
\end{Lem}

\begin{proof}
Since $a^{-1} < 1 < a$, we see that $a^{-1}g < g < ag$.
If $g < x < ag$, then $1 < xg^{-1} < a$, a contradiction.
Hence there is no element strictly
between $g$ and $ag$.  Similarly there is no element between $a^{-1}g$
and $g$.  In particular, for any integers $n<m$, $a^n < g < a^m$
implies $g \in \langle a \rangle$ and thus $\langle a \rangle$ is
convex.  If there exists $f < g$ with no element strictly between, a
routine calculation shows $gf^{-1}$ is the least positive element for
the order.
\end{proof}

Note that according to our definitions, the trivial group 1 has
exactly one right order, and this order is dense but not discrete.

Situated strictly between the class of right-orderable groups and the
class of orderable groups is the class of locally indicable groups.
Recall that a group $G$ is locally indicable if every finitely
generated non-trivial subgroup of $G$ has an infinite cyclic quotient.
Such groups are right-orderable, as was shown by Burns and Hale
\cite{BurnsHale72}.  On the other hand,
a right-orderable group need not be locally
indicable as was shown by Bergman \cite{Bergman91}.  However
for a large class of groups the class of right-orderable groups
coincides with the class of locally indicable groups.  Further results
on this topic are contained in \cite{Linnell01}, \cite{LMR00}, and
especially \cite{Witte06}.

Our interest in considering locally indicable groups $G$ is due to the
fact that such groups have a series (defined below)
with torsion-free abelian factors
as shown by Brodskii in \cite{Brodskii84}.  They possess a right-order
in which the set of convex subgroups form a series with factors which
are order isomorphic to subgroups of the additive group of reals;
we shall refer to such orders as \emph{lexicographic}.
(Such orders are also called Conrad orders and are characterized by
the condition that if $g>1$ and $h>1$, then there exists a positive
integer $n$ such that $(gh)^n > hg$; see
\cite[\S 7.4]{MuraRhemtulla77} for further details.)
Note that for such an ordering, any
non-trivial element $g \in G$ is positive if and only if the
cosets satisfy
$Cg > C$ in the factor group $D/C$ determined by $g$.

By a \emph{series} for $G$ we mean a set
$\Sigma = \{H_\lambda \mid \lambda \in \Lambda \}$ of subgroups of
$G$, where $\Lambda$ is a totally ordered set of indices, satisfying:
\begin{itemize}
\item
if $\lambda < \mu$ then $H_\lambda \subset H_\mu$,

\item $\{1\}$ and $G$ belong to $\Sigma$,

\item $\Sigma$ is closed under arbitrary unions and intersections,

\item if $\mu$ immediately follows $\lambda$ in $\Lambda$, then
$H_\lambda$ is normal in $H_\mu$ and $H_\mu / H_\lambda$ is called a
\emph{factor} associated to the jump $H_\lambda \rightarrowtail
H_\mu$.
\end{itemize}

In the next section we characterize groups that have discrete orders.
We show that a group $G$ has a discrete order if and only if it is an
orderable group and its center $Z(G)$ contains an isolated infinite
cyclic group.  Recall that a subgroup $H$ of a group $G$ is said to be
isolated if $g \in G$ and $g^n \in H$ for some $n > 0$ implies $g \in
H$.

In Section \ref{Sdro}, we deal with groups possessing discrete
lexicographic right-orders and discrete right-orders.  It will follow,
in particular, that any finitely generated
orderable group has discrete
right-orders and if it has a central order (as is the case for free
groups, pure braid groups and wreath products or free products of such
groups), then it has discrete lexicographic right-orders.  Recall that
an order $<$ on $G$ is called central if for every convex jump $C
\rightarrowtail D$, we have $[D, G] \subseteq C$ where $[D,G]$
denotes the subgroup $\langle d^{-1}g^{-1}dg \mid d\in D,\ g \in
G\rangle$.

The result is of course not true for orderable groups in general.  The
additive group of rational numbers has no discrete right-order.

The final section presents examples of discretely ordered groups which
have nontrivial subgroups (for example, the commutator subgroup) upon
which the restriction of the given order is dense.
We also note that there exist finitely
generated right-orderable groups, e.g.\ the Artin braid groups
$B_n$, $n \ge 5$, that are not locally indicable, yet have a discrete
right-order.

\section{Discrete Orders}

\begin{Thm} \label{Theorem-1}
If $<$ is a discrete order on a
group $G$, then there exists an element $z$ in the center $Z(G)$ such
that $\langle z \rangle$ is convex under $<$ and $1 \rightarrowtail
\langle z \rangle$ is a jump.  Conversely, if $G$ is an orderable
group and $Z(G)$ contains an isolated infinite cyclic group, then
there is a discrete order on $G$.
\end{Thm}

\begin{proof}
Let $<$ be a discrete order on $G$ with $z > 1$ as the minimal
positive element.  Then $g^{-1}zg$ is positive for every $g \in G$.
Moreover, $z < g^{-1}zg$ implies $1 < gzg^{-1} < z$, a contradiction.
Thus $z \in Z(G)$.  Also Lemma~\ref{Ldense} shows
that $\langle z \rangle$ is convex under $<$ and $1 \rightarrowtail
\langle z \rangle$ is a jump.

Conversely, let $\langle z \rangle$ be an isolated subgroup in the
center $Z(G)$ of an orderable group $G$.  Since $G$ is orderable,
so is $G/Z(G)$, see \cite[Theorem 2.2.4]{MuraRhemtulla77}.  Moreover,
$Z(G)/\langle z \rangle$ is orderable since $\langle z \rangle$ is
isolated in $Z(G)$.  Order $\langle z \rangle$ (with $z$ positive),
$Z(G)/\langle z \rangle$ and $G/Z(G)$.  Now order $G$ as follows.  If
$1 \neq g \in G\setminus Z(G)$, then put $g$ in the positive cone if
$gZ(G)$ is positive; if $g \in Z(G) \setminus \langle z \rangle$,
then put $g$ in the positive cone if $g\langle z \rangle$ is
positive; if $g = z^n$, then put $g$ in the positive cone if $0 < n$.
It is routine to verify that this gives a discrete order on $G$ with
$z$ as the minimal positive element.
\end{proof}

\begin{Cor}\label{CZtimesG}
For any orderable group
$G$, the group $\mathbb{Z} \times G$ has a discrete order.
In particular, every orderable group
embeds in a discretely orderable group, whose
order extends the given order.
\end{Cor}

\section{Discrete Right Orders} \label{Sdro}

We begin this section with the following result which is easy to
prove.  It is not required in the proofs of the other results.

\begin{Lem}If $(G, <)$ is a nontrivial right-ordered group such that
the order $<$ is a well order on the set of positive elements of $G$,
then G is infinite cyclic.
\end{Lem}

\begin{Lem}\label{L3.2} If $<$ is a discrete right order on $G$ and
$a$ the least positive element under $<$ then for any element $1<g
\in G$, we have $1<aga^{-1}$ and $1<a^{-1}ga$.
\end{Lem}

\begin{proof}
Since $a \le g$, $1\le ga^{-1}$.
Thus $aga^{-1}$ is a product of two
positive elements and hence positive.  By Lemma \ref{Ldense}, there
is no element of $G$ strictly between $a^{-1}g$ and $g$.  Since $1 <
g$, $1 \le a^{-1}g$ and so $a \le a^{-1}ga$.
\end{proof}

\begin{Def}
Let $<$ be a right order on a group
$G$, $C$ a subgroup of $G$ and $a \in G$.  We shall say
``\emph{conjugation by $a$ preserves order on $C$}"
to mean that $C$ is
normalized by $\langle a \rangle$ and conjugation by $a$ and by
$a^{-1}$ preserves the order on $(C, <)$.
\end{Def}

\begin{Lem}\label{L3.4}
Suppose $<$ is a right order on a group $G$,
$C$ is a subgroup of $G$, $1\ne a \in G$ and $C \cap \langle a
\rangle = 1$.  If conjugation by $a$ preserves order on $C$, then
there is a discrete right order on the subgroup $\langle C, a
\rangle$ with $a$ as the minimal positive element.  Moreover, this
right order and the given right order agree on $C$.  Finally
if $aEa^{-1} = E$ for all convex subgroups $E$ of $C$, then the
convex subgroups of $H$ under this new right order are $\{1\}$ and
$\langle a,E\rangle$, where $E$ is a convex subgroup of $C$.
\end{Lem}

\begin{proof}
Set $H = \langle a, C\rangle$.
An element $g \in H$ has a unique expression as $g = a^nc$ where
$c \in C$ and $n\in \mathbb{Z}$.
Define the set $P$ as follows: $g \in P$ if $1 < c$ or $c
= 1$ and $n > 0$.  Note that $P \cup P^{-1} = H \setminus
\{1\}$ and $P \cap P^{-1} =
\emptyset$.  Moreover if $g = a^nc$ and $h = a^md$ are in $P$, then
their product $gh = a^{n+m} (a^{-m}ca^m)d \in P$ as conjugation by
$a^m$ preserves order on $C$.  Thus the order $\prec$ on $H$ given by
$g \prec h$ if and only if $hg^{-1} \in P$ is a right order on $H$.
Furthermore if $g = a^nc$ and $h = a^md$ are in $H$, then $g\prec h$
if and only if $c<d$ or $c=d$ and $n<m$.  It is now clear that $a$
is the least positive element under this order, and that $<$ and
$\prec$ agree on $C$.  Finally we verify that the convex subgroups
are $\langle a,E\rangle$, where $E$ is a convex subgroup of $C$.

Set $A = \langle a\rangle$, so $H = AC$.
If $K$ is a nontrivial convex subgroup of $H$,
then $a \in K$ and $C\cap K$ is a convex subgroup of $C$, and we have
$(K\cap C)A = K \cap CA = K$.  Thus $K = \langle a,E\rangle$ where $E
= C\cap K$.  On the other hand if $E$ is a convex subgroup of $C$, we
claim that $AE$ is a convex subgroup of $H$.  Suppose
$a^mb \prec a^nc \prec a^pd\in EA$, where $b,d \in E$, $c \in C$,
and $m,n,p \in \mathbb{Z}$.  Then $b\le c \le d$ and hence $c\in E$,
and it follows that $AE$ is a convex subgroup of $(H,\prec)$, as
required.
\end{proof}

\begin{Thm} \label{Tdislex}
Let $(G,<)$ be an ordered group, $1 \neq a \in G$ and
$C \rightarrowtail D$ the convex jump determined by $a$ (thus $a \in D
\setminus C$ and $D/C$ is torsion-free abelian).  If $D/\langle
a,C\rangle$ is torsion free, then there is a discrete right order on
$G$ with $a$ as the minimal positive element.  Moreover, if also
$[a, F] \subseteq E$ for every jump $E \rightarrowtail F$, then there
is a discrete lexicographic right order on $G$ (with $a$ as minimal
positive element).
\end{Thm}

\begin{proof}
Set $H = \langle a,C \rangle$.  The hypothesis of Lemma \ref{L3.4}
applies and we right order $H$ as described there.
Next we order the factor group $D/H$.  This is possible since $D/H$ is
torsion-free abelian.  Define the set $Q \subset G$ as follows.  If $g
\in H$, then $g \in Q$ if $g$ is positive in the order on $H$
described above.  If $g \in D \setminus H$, then put $g$ in $Q$ if
$gH$ is positive in the order on $D/H$ given above.  If $g \in G
\setminus D$, then put $g$ in $Q$ if $g$ is positive in $(G, <)$,
the original order on $G$.

It is routine to verify that $Q \cup Q^{-1} = G
\setminus \{1\}$, $Q \cap Q^{-1} =
\emptyset$ and $QQ \subseteq Q$, thus giving a right order $\prec$ on
$G$ with $Q = \{ g \in G \mid 1 \prec g\}$.

The same right order $\prec$ is lexicographic if $[a, F] \subseteq E$
for every jump $E \rightarrowtail F \subseteq C$.  The convex
subgroups are $\{1\}$, $\langle a \rangle$, and $\langle a, E\rangle$
for every subgroup $E$ convex under the original order $<$.
This follows from Lemma \ref{L3.4}: note that
$E = \langle a, E\rangle$ if $E \geq D$ and $E\cap C$ is a convex
subgroup of $C$.
\end{proof}

\begin{Cor} \label{Cfree}
Nontrivial free groups have discrete lexicographic right orders.
\end{Cor}
\begin{proof}
This follows from the fact that the descending lower central
series terminates in $\{1\}$ and the factors are free abelian groups.
Thus any element may be made to be the least positive element so long
as it is a primitive element of the factor group that is determined
by the element.
\end{proof}

Corollary \ref{Cfree} can be generalized to free partially
commutative groups.  These are described in \cite[\S
1.1]{DuchampKrob92}, and the definition given there does not
require these groups to be finitely generated.  Free partially
commutative groups are known under many other names, in particular
they are also called right-angled Artin groups \cite{Charney07}, at
least for finitely generated groups.

\begin{Cor}
Nontrivial free partially commutative groups have discrete
lexicographic right orders.
\end{Cor}
\begin{proof}
Free partially commutative groups are residually
nilpotent by \cite[Theorem 2.3]{DuchampKrob92}.  Furthermore
\cite[Theorems 1.1, 2.1]{DuchampKrob92} show that the quotients of
the lower central series are free abelian groups.  The result now
follows from Theorem \ref{Tdislex}.
\end{proof}

Note that a non-abelian free group does not have a
discrete order.  This follows from Theorem \ref{Theorem-1}.

Finally, all surface groups (orientable or not) except the Klein
bottle and projective plane are residually torsion-free nilpotent, by
\cite[Theorem 1]{Baumslag68} (we would like to thank Warren Dicks for
this reference).  Thus these surface groups also have lexicographic
discrete right orders.  With the exception of the torus, these groups
have trivial center and therefore do not enjoy discrete orders.

The pure braid groups $P_n$, like free groups and surface groups, are
also residually torsion-free nilpotent.  But, unlike those examples,
the groups $P_n$ do have discrete orders.  The center $Z(P_n)$ is
infinite cyclic, generated by $z =$ the full twist braid (often
denoted $\Delta_n^2$).  Since $\langle z \rangle$ is trivially
isolated in $Z(P_n)$,
the second part of Theorem \ref{Theorem-1} provides a discrete order
with $\Delta_n^2$ as least positive element.  In fact any discrete
order of $P_n$ must have $\Delta_n^2$ (or its inverse) as least
positive element.

\section{Examples}

A group may have a lexicographic right order and not have any
discrete right order even when the factors formed by the convex jumps
are all infinite cyclic.  One example of this is the following.

\begin{Ex}\label{EG}
Let $G = \langle a_i\mid i\in \mathbb{Z}\rangle$
with defining relations $[a_i, a_j] = 1$ if $|i-j|>1$ and
$a_{i+1}a_ia_{i+1}^{-1} = a_i^{-1}$.

Every right order on $G$ is lexicographic with the subgroups
$\langle a_i\mid i < j \rangle$ forming the chain of convex
subgroups, and every right order is determined by the $a_i$ (i.e.\
whether or not $a_i$ is in the positive cone for each $i\in
\mathbb{Z}$).  This construction is just the expansion of the well
known (Klein bottle) group $D =
\langle a, b \rangle$ where $b^{-1}ab = a^{-1}$.  There are
exactly four right orders on $D$, every one
discrete with $a$ or $a^{-1}$ as the minimal positive element.
\end{Ex}

We will next show that any infinite cyclic extension of the group $G$
of this example
has a discrete right order if it is finitely generated.  However we
can have a meta-cyclic extension of $G$ that is finitely generated
and right orderable but without any discrete right order.  These are
given as Proposition \ref{PG<t>} and Example \ref{E<G,u,v>}.
If $x,t$ are elements of a group, then $x^t$ will denote $t^{-1}xt$.

\begin{Prop}\label{PG<t>}
Let $\Gamma = G\langle t\rangle$ be a
finitely generated infinite cyclic extension of the group $G$ given
in Example \ref{EG}.  Then $\Gamma$ has a discrete right order with
$t$ (or $t^{-1}$) as the minimal positive element.
\end{Prop}

\begin{proof}
Every non-trivial element $g \in G$ has unique expression of the
form $g = a_{r_1}^{d_1}\dots a_{r_k}^{d_k}$ where $r_1 < \cdots <
r_k$ and $d_i \neq 0$ for all $i$.
Call $a_{r_k}^{d_k}$ the leading term of $g$
and denote it by $\ell(g)$.  Call $r_k$ the leading suffix of $g$.

Note that $\ell(g^n) = (\ell(g))^n$ for all $n \in \mathbb{Z}
\setminus \{0\}$.  Moreover, if $\ell(g) =
a_{r}^{j}$, $\ell(h) = a_{s}^{k}$ and $r < s$, then $\ell(gh) =
\ell(hg) = \ell (h)$.  Since $(a_{i+1}^t)^{-1}
(a_i^t) (a_{i+1}^t) = (a_i^t)^{-1}$,
we see that the leading suffix of $(a_{i+1})^t$
is greater than that of $(a_i)^t$ by at least one.  Thus also the
leading suffix of $(a_{i+1})^{t^{-1}}$ is greater than that of
$(a_i)^{t^{-1}}$ by at least one, and we deduce that the leading
suffix of $(a_{i+1})^t$ is greater than
that of $(a_i)^t$ by exactly one.  It follows that $i > j$ implies
that the leading suffix of $(a_i)^t$ is greater than that of
$(a_j)^t$ by exactly $i-j$.

Since $\Gamma$ is finitely generated, the leading suffix of $a_i^t$
(or that of $a_i^{t^{-1}}$) is greater than $i$ for at least one
value of $i$ -- otherwise $\langle a_j \mid j\le i\rangle$ is normal
in $\Gamma$ for every $i$, and hence $\Gamma$ can not be finitely
generated.  Thus if the leading suffix of $a_0^t$ is $n$,
then the leading suffix of $a_i^t$ is $i+n$ for every integer $i$,
and we may assume that $n>0$.

We now right order the group $G$ by putting $a_0, a_1, \dots,
a_{n-1}$ in the positive cone $P$.  Next,
for all $n \leq r <2n$, we
put $a_r \in P$ if the exponent of $\ell(a_{r-n}^t)$ is positive and
$a_r^{-1} \in P$ otherwise.  Next put $a_{r+n}$ or $a_{r+n}^{-1}$ in
$P$ depending on whether the exponent of $\ell((\ell(a_{r-n})^t)^t)$
is positive or negative.  Continue this process.  For every integer $i
\geq 0$ we have determined whether $a_i$ or its inverse is in $P$.
Next, for $0 > r \geq -n$ put $a_r \in P$ if the exponent of
$\ell(a_{r+n}^{t^{-1}})$ is positive.  Put $a_r^{-1} \in P$ otherwise.
Continue this process.
This takes care of all $a_i$ for $i \in \mathbb{Z}$.
Note that the above order on $G$ is $\langle t \rangle$ invariant.
Hence by Lemma \ref{L3.4}, $\Gamma$ has a right discrete order with
$t$ as the minimal positive element.
\end{proof}

\begin{Ex}\label{E<G,u,v>}
Let $G$ be the group in Example \ref{EG}.  Consider the
map $\phi\colon \{a_i\mid i\in \mathbb{Z}\} \to \{a_i^{-1}\mid
i\in \mathbb{Z}\}$ given by $\phi(a_i) = a_i^{-1}$.  Then $\phi$
extends uniquely to an automorphism of $G$ that inverts every $a_i$.
Let $\langle G, u\rangle$ be the infinite cyclic extension of $G$ by
$\langle u \rangle$ where $u^{-1}a_iu = a_i^{-1}$ for all $i \in
\mathbb{Z}$.  Note that $\langle G, u\rangle$ is right orderable
because it is an infinite cyclic extension of
the right orderable group $G$.  However, it has no discrete right
order.  This can be seen as follows.  Suppose $g \in G$ and $c :=
gu^j$ is a minimal positive element under some right order $<$ on
$\langle G,u\rangle$.  Since $G$ has no discrete right order,
$j \neq 0$ and it must
be even otherwise $c^{-1}a_ic = a_i^{-1}$ for some $i$, contradicting
Lemma \ref{L3.2}.

Suppose $1 < u$.  Then $a_i^{r} < u$ for every $i, r \in \mathbb{Z}$,
for if $1 < u < a_i^{r}$, then $ua_i^{-r} < 1$.  Hence $a_i^{r} =
ua_i^{-r}u^{-1} < 1$, a contradiction.  Thus $h < u$ for all $h \in
G$.  Since $1 < c = gu^j = u^jg$, we see that $j$ is positive, and
then we have $1 < u^{j-1}, ug$, which contradicts the hypothesis
that $c$ is the minimal positive element.
The argument is similar if $u < 1$.  We note in particular that if
$1<u$, then $h<u$ for every $h \in G$.

Next extend the group $\langle G, u\rangle$ by
the infinite cyclic group
$\langle v \rangle$ to get the group $J =\langle G,u,v\rangle$ where
the action of $v$ under conjugation is as follows: $v^{-1}a_iv =
a_{i+1}$, the shift automorphism, and $v^{-1}uv = u^{-1}$.  Note that
$\langle G, v\rangle = \langle a_0, v\rangle$, $J = \langle a_0, u,
v\rangle$ and $J$ is right orderable.

We now show that $J$ has no
discrete right order.  Suppose that $c := gu^jv^k$ is a minimal
positive element under a right order $<$ on $J$.
Then $k \neq 0$, since otherwise
the restriction of the right order $<$ to $\langle G, u\rangle$ would
be discrete with $gu^j$ as the minimal positive element.  We have
seen that this is not possible.  Next note that $k$ must be even,
otherwise assume without loss of generality that $u>1$.  Then
$cuc^{-1} = gu^{-1}g^{-1}< 1$, which contradicts Lemma \ref{L3.2}.
Since $vuv^{-1} = u^{-1}$, we see that $1 < v$
implies $x < v$ for every $x \in \langle G, u\rangle$, in particular
$k>0$ and $gu^jv^{k-1} > 1$.  This contradicts the hypothesis that
$c$ is the minimal positive element.  Similarly we cannot have $1>v$,
which finishes the verification that $J$ has no discrete right order.
\end{Ex}

It is possible for a discretely (right-) ordered group to have a
subgroup on which the same order is dense.  Indeed, by Corollary
\ref{CZtimesG}, any densely ordered group is a subgroup of a
discretely ordered group, whose order extends the given order.
Following is a ``natural'' example of this phenomenon for
right-ordered groups.

\begin{Ex}
The Artin braid groups $B_n$ have a discrete
right-order, which becomes dense when restricted to the commutator
subgroup.  For each integer $n \ge 2$, $B_n$ is the group generated
by $\sigma_1, \sigma_2, \dots, \sigma_{n-1}$, subject to the relations
\[
\sigma_i\sigma_j = \sigma_j\sigma_i \text{ if } |i-j| >1,\quad
\sigma_i\sigma_j\sigma_i=\sigma_j\sigma_i\sigma_j
\text{ if } |i-j| =1.
\]

It was shown by Dehornoy (see \cite{Dehornoy94} and \cite{DDRW02})
that each $B_n$ is right-orderable (but not orderable, for $n>2$).
The positive cone consists of all elements expressible as a word in
the $\sigma_i$ such that the generator
with the lowest subscript occurs
with only positive exponents.  This right-order is discrete, with
smallest positive element $\sigma_{n-1}$.  On the other hand, it is
shown in \cite{ClayRolfsen07} that the Dehornoy order, when
restricted to the commutator subgroup $B'_n = [B_n,B_n]$, is a dense
order for $n \ge 3$.  For $n =3$, $B'_n$ is free (on two
generators).  For $n \ge 5$, $B'_n$ is finitely-generated and
perfect (see \cite{GorinLin69}), so $B_n$ is an example of a
non-locally indicable discretely right-orderable group for $n \ge
5$.

Now consider the braid group $B_3$ with its two generators $\sigma_1$
and $\sigma_2$ and let $H$ be the subgroup generated by
$\sigma_1^2$ and $\sigma_2^2$.  Crisp and Paris \cite{CrispParis01}
showed that $H$ is a free group with free basis $\sigma_1^2$ and
$\sigma_2^2$.  The Dehornoy
order restricted to this subgroup has the least positive element
$\sigma_2^2$.  This gives an alternative construction of discrete
right-orders on a free group.
\end{Ex}

Bergman \cite{Bergman91} published the first examples of
groups which are right-orderable and not locally indicable; some of
his examples are finitely generated and perfect.  We shall argue that
they can be given a discrete right-order.

If $G$ is a group acting on a set and $x_1,\dots,x_n$ are elements of
the set, then $\Stab_G(x_1,\dots,x_n)$ will denote the pointwise
stabilizer of $\{x_1,\dots,x_n\}$ in $G$, namely $\{g \in G \mid gx_i
= x_i$ for all $i\}$.  Also $I$ will denote the identity matrix of
$\SL_2(\mathbb{R})$.
We have an action of $\SL_2(\mathbb{R})$ on the one
point compactification $\overline{\mathbb{R}} =\mathbb{R} \cup
\{\infty\} \cong S^1$, the circle, given by the rule
\begin{equation} \label{Eaction}
\begin{pmatrix} a&b\\ c&d \end{pmatrix}(x) = \frac{ax+b}{cx+d}.
\end{equation}
This induces a faithful action of $\PSL_2(\mathbb{R})$ on
$\overline{\mathbb{R}}$.
Let $H$ be a finitely generated subgroup of $\SL_2(\mathbb{R})$
containing the center $\{\pm I\}$ of $\SL_2(\mathbb{R})$,
and let $\bar{H}$ denote its image in $\PSL_2(\mathbb{R})$.
Since $\mathbb{R}$ is the universal covering space of $S^1$,
we can lift the action of $\bar{H}$ on $S^1$ to an action of a group
$G$ on $\mathbb{R}$ by orientation preserving homeomorphisms.
In this situation, $G$
will have a central subgroup $Z \cong \mathbb{Z}$ such that $G/Z
\cong \bar{H}$ and $Z$ acts fixed point free on $\mathbb{R}$.
Also if $\pi\colon \mathbb{R} \to S^1$ is the associated covering
map and $p \in \mathbb{R}$, then $\Stab_{\bar{H}}(\pi p) =
Z\Stab_G(p)/Z \cong \Stab_G(p)$.

\begin{Prop} \label{PBergmandiscrete}
Let $H$ be a finitely generated subgroup of $\SL_2(\mathbb{R})$
with $-I \in H$ and
let $G$ be its lift to orientation preserving homeomorphisms of
$\mathbb{R}$ (as described above).  Suppose $H$ contains a diagonal
matrix other than $\pm I$.  Then $G$ has a discrete right order.
\end{Prop}

To prove this, we will need the following; for a proof, see
\cite[Lemma 2.2]{Linnell06}.
\begin{Lem} \label{Lextend}
Let $G$ be a right ordered group, let $H$ be a convex subgroup
of $G$ and let $<$ be any right order on $H$.  Then there exists a
right order on $G$ whose restriction to $H$ is $<$, and $H$ is still
a convex subgroup under this new right order.
\end{Lem}
\begin{proof}[Proof of Proposition \ref{PBergmandiscrete}]
Let us examine
$\Stab_H(0)$ and $\Stab_H(0,\infty)$ with the action given by
\eqref{Eaction}.  The former is the lower triangular
matrices, i.e.~the matrices above with $b=0$; we shall denote by $L$
those lower triangular matrices which lie in $H$.  The latter is
given by the diagonal matrices; we shall denote by $D$ those diagonal
matrices which lie in $H$.

Thus we have an action of $H$ on $S^1$ and two points $p_1,p_2 \in
S^1$ such that $\Stab_H(p_1) = L$ and $\Stab_H(p_1,p_2) = D$.
Let $p_3 \in S^1$ be distinct from $p_1,p_2$.  Then
$\Stab_H(p_1,p_2,p_3) = \{\pm I\}$.

Now we can lift the action of $H$ on $S^1$ to an action of
$G$ on $\mathbb{R}$ by orientation preserving homeomorphisms; $G$
will have a central subgroup $Z \cong \mathbb{Z}$ such that $G/Z
\cong H/\{\pm I\}$.  For $i = 1,2,3$,
let $q_i \in \mathbb{R}$ be a lift of
$p_i$.  Then $Q:= \Stab_G(q_1,q_2) \cong D/\{\pm I\}$.  We now
define a right order on $G$ in the usual way when we have a group
acting on $\mathbb{R}$.  The positive cone is the set of all $g \in
G$ such that $g(q_i) > q_i$ for the smallest $i$ such that $g(q_i)
\ne q_i$.  This right order will have the property that $Q$ is a
(smallest nontrivial) convex subgroup of $G$.  If $Q \cong
\mathbb{Z}$, then it would follow that the above defined right
order will be discrete, but this is not true in general.
However since $H$ is finitely generated, $H
\subseteq \SL_2(R)$ for some finitely generated subring $R$ of
$\mathbb{R}$.  By \cite[Th\'eor\`eme 1]{Samuel66}, the group of units
of a finitely generated integral domain is finitely generated, hence
$D$ is also finitely generated.  We deduce that
$D/\{\pm I\}$ is a finitely generated free abelian group.
Thus $Q$ is also a finitely generated
free abelian group and hence has a discrete right order by Corollary
\ref{CZtimesG}.  The result now follows from Lemma \ref{Lextend}.
\end{proof}

\bibliographystyle{plain}

\end{document}